\documentclass[letterpaper, 10 pt, conference, twocolumn]{ieeeconf}  % Comment this line out

\IEEEoverridecommandlockouts
\usepackage{graphics}
\usepackage{epsfig}
\usepackage{amsmath}
\usepackage{amssymb}
\usepackage{bbm}
\usepackage{eufrak}
\usepackage[noadjust]{cite}

\newtheorem{theorem}{Theorem}

\newtheorem{remark}{Remark}

\usepackage{comment}
\usepackage{mathtools}
\usepackage{soul}
\usepackage{tikz}
\usepackage{color}
\usepackage{xcolor}
\usepackage{tcolorbox}
\usepackage{float}
\usepackage{subfig}
\usepackage{algorithm}
\usepackage[noend]{algpseudocode}
\usepackage[noabbrev]{cleveref}
\usepackage{scalerel,amssymb}

\makeatletter
\def\BState{\State\hskip-\ALG@thistlm}
\makeatother

\graphicspath{{Images/}}
\title{\LARGE \bf Distributed $Q$-Learning for Dynamically Decoupled Systems}

\author{Siavash Alemzadeh and Mehran Mesbahi
	\thanks{The authors are with the William E. Boeing Department of	Aeronautics and Astronautics, University of Washington, Seattle WA, USA {\tt\small \{alems,mesbahi\}@uw.edu}.
	The research of the authors has been supported by NSF grant SES-1541025 and ARO Grant W911NF-13-1-0340.}
}

\begin{document}
	
	\maketitle
	\thispagestyle{empty}
	\pagestyle{empty}
	
	%=================================
	
	\begin{abstract}
		Control of large-scale networked systems often necessitates the availability of complex models for the interactions amongst the agents.
		However, in many applications building accurate models of these interactions might be prohibitive due to the curse of dimensionality or their inherent complexity.
		In the meantime, data-guided control methods can circumvent model complexity by directly synthesizing the controller from the observed data.
		In this paper, we propose a distributed $Q$-learning algorithm to design a feedback mechanism given an underlying graph structure parameterizing the agents' communication.
		We assume that the distributed nature of the system arises from a common cost and show that for the particular case of identical dynamically decoupled systems, the learned controller converges to the optimal Linear Quadratic Regulator controller for each subsystem.
		We provide a convergence analysis and verify the result with an example.
		
	\end{abstract}
	
	%=================================
	
	\section{INTRODUCTION}
	\label{sec:intro}
	
	Distributed control has undergone an unprecedented growth during the past few years mainly due to the complexity in modeling and analysis of large scale systems.
	In such scenarios, high-dimensional collective tasks conducted by members of a team are formed from local decisions of each member leading towards the global system-level final decision.
%	One major inspiration of the field is due to high-dimensional collective tasks conducted by the members of a team pursuing a common goal.
%	These members can be individuals in a social community, mobile robots in a swarm, or birds flying in flocks.
%	Although these units form a centralized body, the ultimate convergent pattern is essentially formed from local decisions of each member leading towards the global system-level final decision.
	Accordingly, the main focus in distributed control design is finding the closest-to-optimal control mechanism for a large-scale system, making use of the structure in information exchange and decision-making.
%	This sparsity comes from lack of information exchange between agents or groups. 
	Indeed, such an approach has found broad applications in such areas as robotic swarms \cite{bullo2009distributed}, structured robust learning \cite{sedghi2019robust}, and social networks \cite{alemzadeh2017controllability}.
	
	Decentralized control of large-scale systems is not a new research area.
	The roots of the field traces back to the socioeconomics literature of 1970's \cite{mcfadden1969controllability}; an early work in the control literature is~\cite{wang1973stabilization}.
	The inspiration of these types of works is that the presupposition of centrality fails to hold due to the lack of either central intelligence or computational capability~\cite{sandell1978survey}.
	This line of work was followed by the pioneering work \cite{corfmat1976decentralized}, where stability conditions for multi-channel linear systems were derived.
	Fast forward a few decades, the stability of networks was studied in \cite{fax2004information}, where sufficient graph-theoretic conditions were provided for stability of formations comprised of identical vehicles.
	Graph-based structured controller design was further examined in works such as \cite{borrelli2008distributed, massioni2009distributed}.
	The topic was also studied from a spatially distributed control viewpoint~\cite{bamieh2002distributed} or a layered control design approach~\cite{alemzadeh2018influence}.
%	Distributed control is also studied from other viewpoints in the literature.
%	For example, the work~\cite{bamieh2002distributed, motee2008optimal} considers spatially distributed control and a layered control design approach has been explored in works such as~\cite{chapman2017data, alemzadeh2018influence}.
	However, all of these studies are based on the knowledge of the underlying dynamics; as the system grows in scale, modeling becomes prohibitively difficult and uncertain due to complexities or potential perturbations in high dimensions.\footnote{For example, computational complexity of order $\mathcal{O}(n^3)$ for solving the Algebraic Riccati Equation (ARE) is not desirable for large-scale systems.}
	This motivates a data-guided approach to evade the difficulties of model-based distributed control.
	
	In this work, we focus on a model-free distributed control design using $Q$-learning \cite{watkins1992q}.
	The approach has been used to find the optimal Linear Quadratic Regulator (LQR) feedback controller online for a single system \cite{bradtke1994incremental}.
	In that sense, $Q$-learning can be thought of as an adaptive optimal control design method \cite{lewis2012reinforcement}.
	More recent works on distributed adaptive systems can be found in \cite{gao2017data,nguyen2017selectively} to name a few.
	Other related works are~\cite{kar2013cal} that considers network effects within $Q$-functions, and \cite{nguyen2017selectively} that introduces a decentralized $Q$-learning approach for a general framework but not necessarily on an underlying graph.
	Moreover, there has emerged several recent works on data-driven control using finite sample opposed to the adaptive control design in the limit \cite{fattahisample,dean2017sample,alaeddini2018linear}.
	
	Our contribution is mainly built upon the work of Bradtke~\cite{bradtke1994incremental}, in that we provide a distributed policy iteration to find a collective controller.
	We assume that the distributed nature of the problem comes from the interconnection of identical dynamically decoupled systems that work together for a common system-level goal in a network.
	In fact, we assume that the underlying interaction graph is only reflected in the performance index of the corresponding LQR problem.
%	\textcolor{red}{It is shown that this type of modeling can describe applications such as formation flight and monitoring cameras \cite{borrelli2005hybrid}.
%	The algorithm introduced in this paper is mostly comparable to the one introduced in \cite{bradtke1994incremental}.
%	However, the works have fundamental differences in the way the interconnections of the systems work.
%	Moreover, the section provides no theory about the convergence on the distributed policy iteration algorithm.
%	We are inspired by \cite{bradtke1994incremental} regarding the huge reduction in the computational complexity which will be restated in the analysis.}
	We provide a graph theoretical framework under which each agent synthesizes an estimate of the optimal LQR controller.
	%
%	We also provide a proof of convergence and a computational comparison between the proposed method and the centralized scenario.
%	Furthermore, an illustrative example is provided to validate the algorithm.
	
	The rest of the paper is organized as follows: In \S\ref{sec:math} we provide a quick overview of
	mathematical tools that are used in the paper.
	In \S\ref{sec:p-setup} we introduce the problem setup. In \S\ref{sec:DQL} the distributed setup and the main algorithm is provided along with its proof of convergence.
	The section concludes with a discussion on the computational savings due to the adoption of the distributed algorithm.
	An example is provided in \S\ref{sec:example} to validate our theoretical results.
	Concluding remarks and future directions are discussed in \S\ref{sec:conclusion}.
	
	%=================================
	
	\section{Mathematical Preliminaries}
	\label{sec:math}
	
	We denote by $\mathbb{R}$ the set of real numbers.
	A column vector with $n$ elements is referred to as $v\in\mathbb{R}^n$, where $v_i$ represents the $i$th element in $v$.
	The matrix $M\in \mathbb{R}^{p\times q}$ contains $p$ rows and $q$ columns with $[M]_{ij}$ denoting the element in the $i$th row and $j$th column of $M$.
	The square matrix $N\in\mathbb{R}^{n\times n}$ is \emph{symmetric} if $N^{\top}=N$, where $N^\top$ denotes the \emph{transpose} of the matrix $N$.
	The operator \emph{diag(.)} makes a square diagonal matrix out of the elements of its argument.
	The $n\times n$ zero matrix is denoted by $\textbf{0}_n$ and $\textsc{I}_n = \text{diag}(1,1,\dots,1)$, is the \emph{identity matrix}.
	We write $N\succ0$ ($\succeq0$) when $N$ is a positive-(semi)definite matrix, i.e., $x^{\top}Nx>0$ ($\geq0$) for all $x\neq0$.
	To simplify the vector notation, we use semicolon (;) to concatenate column vectors, hence $[v^{\top}\ w^{\top}]^{\top} = [v;w]$.
	We call the pair $(A,B)$ \emph{controllable}, if and only if the controllability matrix $\mathcal{C}=[B\quad AB\ \dots\ A^{n-1}B]$ has full-rank, where $n$ is the size of the system.
	$A\otimes B\in\mathbb{R}^{p_1q_1\times p_2q_2}$ refers to the \emph{Kronecker product} of $A\in\mathbb{R}^{p_1\times q_1}$ and $B\in\mathbb{R}^{p_2\times q_2}$, and $\textsc{I}\otimes T$ gives a block diagonal square matrix, with $T$ on each diagonal block.
	A graph is characterized by $\mathcal{G}=(\mathcal{V},\mathcal{E})$ where $\mathcal{V}$ is the set of nodes and $\mathcal{E}\subseteq \mathcal{V}\times\mathcal{V}$ denotes the set of edges.
	An edge exists from node $i$ to $j$ if $(i,j)\in \mathcal{E}$;
this is also specified by writing $j\in\mathcal{N}_i$, where $\mathcal{N}_i$ is the set of neighbors of node $i$.
	Finally, $\mathcal{G}$ can be represented by various matrices, in particular, by its \emph{graph Laplacian} denoted by $\mathcal{L}$.
	
	%=================================
	
	\section{Problem Setup}
	\label{sec:p-setup}
	
	Herein, we provide the basic formulation and problem setup.
	First, we introduce the distributed LQR control on a given graph.
	This is mainly related to~\cite{borrelli2008distributed}, where the network contains identical dynamic agents,
	yet decoupled from other agents' dynamics.
	The only coupling between these agents is through a common network-level objective function.
	Then we introduce the basic setup of $Q$-learning for linear dynamical systems and extend the formulation for the distributed setup.
	As we shall see later, the distributed nature will be simplified into an additional interaction term in the output of a linear Recursive Least Squares (RLS) algorithm.
	
	\subsection{Distributed LQR Problem}
	\label{subsec:disLQR}
	
	Assume that the system contains $N$ agents that form a graph $\mathcal{G}$	with each node of the graph indicating a linear time-invariant dynamical system corresponding to that agent as,
	\begin{align*}
		x_{t+1}^{(i)} = Ax_t^{(i)} + Bu_t^{(i)},
		\qquad i = 1,2,\dots,N,
	\end{align*}
	where $x_t^{(i)}$ is the state of agent $i$ at time-step $t$, where $x_t^{(i)}\in\mathbb{R}^n$, $u_t^{(i)}\in\mathbb{R}^m$, $A\in\mathbb{R}^{n\times n}$, and $B\in\mathbb{R}^{n\times m}$.
	The assumption that all agents have identical system matrices $A$ and $B$ is relevant in many applications such as formation flight, homogenous mobile robots, and power grids consisting of identical generators.
	These dynamics can be integrated into a compact form as, $\tilde{\textbf{x}}_{t+1} = \tilde{\textbf{A}} \tilde{\textbf{x}}_{t} + \tilde{\textbf{B}} \tilde{\textbf{u}}_{t},$
%	\begin{align*}
%		\tilde{\textbf{x}}_{t+1} = \tilde{\textbf{A}} \tilde{\textbf{x}}_{t} + \tilde{\textbf{B}} \tilde{\textbf{u}}_{t},
%	\end{align*}
	where $\tilde{\textbf{x}}\in\mathbb{R}^{Nn}$ and $\tilde{\textbf{u}}\in\mathbb{R}^{Nm}$ are formed by concatenation of all states and inputs into one vector with $\tilde{\textbf{A}} = \textsc{I}_{N} \otimes A\in\mathbb{R}^{Nn\times Nn}$ and $\tilde{\textbf{B}} = \textsc{I}_{N} \otimes B\in\mathbb{R}^{Nn\times Nm}$.
%	\begin{gather*}
%		\tilde{\textbf{A}} = \textsc{I}_{N} \otimes A, \quad
%		\tilde{\textbf{B}} = \textsc{I}_{N} \otimes B,
%	\end{gather*}
%	with $\tilde{\textbf{A}}\in\mathbb{R}^{Nn\times Nn}$ and $\tilde{\textbf{B}}\in\mathbb{R}^{Nn\times Nm}$.
	The graph structure is reflected in the cost function of the associated LQR problem by the following definition,
	\begin{equation}
		\begin{aligned}
			\textbf{J} &= \sum_{i=1}^{N} \big[ \underbrace{x_t^{(i)\top} Q_i x_t^{(i)} + u_t^{(i)\top} R_i u_t^{(i)} }_{\text{System $i$ Performance}} \big] \\
			& \hspace{8mm}+\sum_{i=1}^{N} \sum_{j\neq i}^{N} \big[ \underbrace{\  \big( x_t^{(i)}-x_t^{(j)} \big)^{\top} Q_{ij} \big( x_t^{(i)} - x_t^{(j)} \big) }_{\text{Systems $i$ and $j$ Coupling}}\big],
		\end{aligned}
		\label{eq:cost}
	\end{equation}
	where the first term indicates the intra-systems cost while the second denotes the inter-system coupling.
	We make the simplifying assumption,
	\begin{align*}
		Q_i = \bar{Q},\quad R_i = \bar{R},\quad Q_{ij} = \bar{Q}_{\{ j\in\mathcal{N}_i \}} = \begin{cases}
			\textbf{0}_n \hspace{5mm} j\not\in\mathcal{N}_i \\
			\bar{Q} \hspace{6mm} j\in\mathcal{N}_i
		\end{cases}
	\end{align*}
	where $\bar{Q}\succeq 0$ and $\bar{R}\succ 0$.
	The cost function can also be written in compact form as $\textbf{J} = \tilde{\textbf{x}}_t^{\top} \tilde{\textbf{Q}} \tilde{\textbf{x}}_t + \tilde{\textbf{u}}_t^{\top} \tilde{\textbf{R}} \tilde{\textbf{u}}_t$
%	\begin{align}
%		\textbf{J} = \tilde{\textbf{x}}_t^{\top} \tilde{\textbf{Q}} \tilde{\textbf{x}}_t + \tilde{\textbf{u}}_t^{\top} \tilde{\textbf{R}} \tilde{\textbf{u}}_t,
%		\label{eq:compactCost}
%	\end{align}
	where, $\tilde{\textbf{Q}} = (\mathcal{L}+\textsc{I}_n) \otimes \bar{Q} \succeq 0$ and $\tilde{\textbf{R}} = \textsc{I}_n \otimes \bar{R} \succ 0$.
%	\begin{align*}
%		\tilde{\textbf{Q}} = (\mathcal{L}+\textsc{I}_n) \otimes \bar{Q}, \qquad \tilde{\textbf{R}} = \textsc{I}_n \otimes \bar{R}.
%	\end{align*}	
	
	Solution of the LQR problem in such systems is studied for a particular $\tilde{\textbf{Q}}$ resulting in a structured controller~\cite{borrelli2008distributed}.
	Suboptimal solutions to the controller design consistent with the graph structure has also been proposed.
	Nevertheless, in many real-world applications there is no \emph{a priori} knowledge of the system's model due to either complexities or model uncertainties~\cite{alemzadeh2018influence}.
	We introduce a model-free approach while considering the optimality criteria for each subsystem.
	We will show that for an interconnected system with identical dynamically decoupled agents as discussed above, $Q$-learning leads to each subsystem running their respective local LQR optimal controller independent of other agents in the network.
	This phenomenon is shown to hold asymptotically after each agent collects enough data.
	
	\vspace{2mm}
	
	\begin{remark}
%		While it is shown that defining the LQR cost in the form of \eqref{eq:cost} can have applications such as formation flight and monitoring cameras \cite{borrelli2005hybrid}, we point out the cost structure can further be is still rudimentary as the coupling term goes to zero trivially because of the existence of the first term in the cost function.
		We note that the global cost in \cref{eq:cost} induces a structured way of steering the states of the agents to the origin through an auxiliary consensus term.
		Our future work will consider further realizations of the global/local cost structure in the LQR setup--that might not be completely aligned with each other.
		%
		%Changing the cost into a more sensible scenario is addressed as a future work.
	\end{remark}
	
	\subsection{Centralized $Q$-Learning}
	\label{subsec:disQfun}
	
	To make the paper self-contained, we refer to some basics of $Q$-learning and its connections to LQR feedback control design.
	$Q$-learning describes a methodology where an agent aims to optimize the value of a sum of reward functions from observing the results of its own actions.
	This value is reformulated by the $Q$-function which is defined for a single agent as,
	\begin{equation}
		\begin{aligned}
			Q(x_t,u_t) = R(x_t,u_t) + \gamma Q(x_{t+1},u_{t+1}),
		\end{aligned}
		\label{eq:Qdefine}
	\end{equation}
	where $Q(x_t,u_t)=x_t^{\top}Px_t$ is the \emph{state-action} $Q$-function, $P$ is the \emph{cost-to-go} matrix,\footnote{which is also the solution to the discrete-time ARE in LQR.} and $R(x_t,u_t)=x_t^{\top}\bar{Q} x_t + u_t^{\top}\bar{R} u_t$ is the \emph{one-step reward} with symmetric constant matrices $\bar{Q}\succeq 0$ and $\bar{R}\succ 0$.
	\Cref{eq:Qdefine} is the simplified form of the well-known Bellman equation for the deterministic case of LQR.
	Also, the control actions come from a set of optimal policies that assume the form of a feedback law $u_{t}=-Kx_{t}$ in the LQR framework.
Simplification of \eqref{eq:Qdefine} results in,
	\begin{equation}
		\begin{aligned}
			Q(x_t,u_t) = z_t^{\top} H z_t,
		\end{aligned}
		\label{eq:Qfun}
	\end{equation}
	where $z_t = [x_t;u_t]$ and $H$ is a block matrix defined as,
	\begin{align*}
		H = \begin{bmatrix}
			H_{11} & H_{12} \\
			H_{21} & H_{22}
		\end{bmatrix}
		= \begin{bmatrix}
			\bar{Q} + \gamma A^{\top}PA & \gamma A^{\top} P B \\
			\gamma B^{\top} P A & \bar{R} + \gamma B^{\top} P B 
		\end{bmatrix}.
	\end{align*}
	Then the idea is to learn the parameters in $H$ through observations $z_t$ and update the estimate of the controller as,
	\begin{align*}
		K_{\text{new}} = -H_{22}^{-1} H_{21} = -\gamma (R + \gamma B^{\top} P B)^{-1} (B^{\top} P A),
	\end{align*}
	which can also be obtained by setting $\partial Q/\partial u_t = 0$.
	The adaptive nature of the algorithm is originated from a linear RLS step to learn the parameters of $H$ in real-time.
	Hence, we pursue \cite{bradtke1994incremental} to form a linear parameterization of \eqref{eq:Qfun} as,
	\begin{align}
		Q(x_t,u_t) = z_t^{\top} H z_t = \bar{z}^{\top}_t \theta_H,
		\label{eq:quad}
	\end{align}
	where $\bar{z}_t,\theta_H\in\mathbb{R}^{(n+m)(n+m+1)/2}$ are quadratic basis of the elements in $z_t$ and vector of upper right triangle of symmetric $H$ in the correct order, respectively.
	With these definitions,
	\begin{equation}
		\begin{aligned}
			R(x_t,u_t) = r_t &= Q(x_t,u_t) - \gamma Q(x_{t+1},u_{t+1}) \\
			&= z_t^{\top} H z_t - \gamma z_{t+1}^{\top} H z_{t+1} = \phi_t^{\top} \theta_H,
		\end{aligned}
		\label{eq:reward}
	\end{equation}
	where $\phi_t = \bar{z}_t - \gamma\bar{z}_{t+1}$.
	Therefore, assuming that we know $R(x_t,u_t)$ and $\phi_t$, RLS can be employed to find an estimate of $\theta_H$.
	According to \cite{goodwin2014adaptive}, this recursive algorithm converges in the limit if $\phi_t$ is persistently excited (PE), i.e.,
	\begin{align}
		\alpha\textsc{I} \leq \frac{1}{M} \sum_{i=1}^M \phi_{t-i}\phi_{t-i}^{\top} \leq \beta\textsc{I} \qquad \forall\  t,M\geq M_0,
		\label{eq:PE}
	\end{align}
	for some positive parameters $M_0$, $\alpha$, and $\beta$.	
	Following the convergence of $\theta_H$, then $H$ is obtained using \eqref{eq:quad}.	
	
	%=================================
	
	\section{Distributed $Q$-learning}
	\label{sec:DQL}
	
	%=================================
	
	\subsection{Distributed $Q$-function}
	\label{sec:DQfun}
	
	We now switch to a multiagent setup, where several autonomous agents try to minimize their own discounted reward based on a global cost and single-agent control is not applicable since there exist multiple decision-makers.
	In this section, we extend the $Q$-learning setup based on the distributed control framework defined in \cref{subsec:disLQR}.
	To this end, we assume that each agent enjoys its own $Q$-function whose reward is a function of the state of the agent as well as the state of its neighbors.
	For agent $i$ we define,
	\begin{equation}
		\begin{aligned}
			Q^{(i)}(x_t^{(i)},u_t^{(i)}) &= R^{(i)}(x_t^{(i)},u_t^{(i)}) + \gamma Q^{(i)}(x_{t+1}^{(i)},u_{t+1}^{(i)}) \\
			&=  y_t^{(i)\top} \mathcal{Q}^{(i)} y_t^{(i)} + \gamma x_{t+1}^{(i)\top} P^{(i)} x_{t+1}^{(i)},
		\end{aligned}
		\label{eq:i}
	\end{equation}
	where $y_t^{(i)} = [x_t^{(i)} ;\ u_t^{(i)} ;\ x_t^{(j_1)} ;\ \dots\ ;\ x_t^{(j_{d_i})}]$, $d_i$ is the degree of agent $i$, $j_{d_k}\in\mathcal{N}_i$ for $k=1,\dots,i$, and $\mathcal{Q}^{(i)}$ is defined as, \footnotesize
	\begin{equation}
		\begin{aligned}
			\mathcal{Q}^{(i)} = \begin{bmatrix}
				(d_i+1)\bar{Q}   & \textbf{0} & -\bar{Q}  & \dots & -\bar{Q} \\
				\textbf{0} & \bar{R} & \textbf{0} & \dots & \textbf{0} \\
				-\bar{Q} & \textbf{0}& \bar{Q} & \dots & \textbf{0} \\
				\vdots   &  \vdots   &   \vdots   & \ddots& \vdots \\
				-\bar{Q} & \textbf{0} & \textbf{0} & \dots& \bar{Q}
			\end{bmatrix}\in\mathbb{R}^{(d_i+2)n\times (d_i+2)n}.
		\end{aligned}
		\label{eq:bigQ}
	\end{equation} \normalsize
	The structure of $\mathcal{Q}^{(i)}$ is resulting from \cref{eq:cost} and implies the new definition of reward function for multiple agents in the system.
%	\footnotesize
%	\begin{align*}
%		y_t^{(i)\top} \mathcal{Q}^{(i)} y_t^{(i)} = x_t^{(i)\top}\bar{Q}x_t^{(i)} + \sum_{k=1}^{d_i} \big( x_t^{(i)}-x_t^{(j_k)} \big)^{\top} \bar{Q} \big( x_t^{(i)}-x_t^{(j_k)} \big)
%	\end{align*} \normalsize
	Note that \cref{eq:i,eq:bigQ} make two implicit assumptions: ($i$) there is no control coupling amongst agents and, ($ii$) each agent has only access to the reward form	the coupling between its own state and the states of neighbors.
	This motivates the existence of zero blocks in \eqref{eq:bigQ}.
	Similar to \eqref{eq:Qfun}, \cref{eq:i} can also be re-arranged into,
	\begin{equation*}
		\begin{aligned}
			Q^{(i)}(x_i,u_i) = y_t^{(i)\top} H^{(i)} y_t^{(i)},
		\end{aligned}
	\end{equation*}
	where, \scriptsize
	\begin{equation*}
		\hspace{-0.8mm}
		\begin{aligned}
			H^{(i)} = \begin{bmatrix}
			(d_i+1)\bar{Q}+\gamma A^{\top}P^{(i)}A & \gamma A^{\top} P^{(i)} B & -\bar{Q} & \dots & -\bar{Q} \\
			\gamma B^{\top} P^{(i)} A & \bar{R} + \gamma B^{\top} P^{(i)} B & \textbf{0} & \dots & \textbf{0} \\
			-\bar{Q} & \textbf{0} & \bar{Q} & \dots & \textbf{0} \\
			\vdots & \vdots & \vdots& \ddots & \vdots \\
			-\bar{Q} & \textbf{0} & \textbf{0} & \dots & \bar{Q}
			\end{bmatrix}
		\end{aligned}
	\end{equation*}
	\normalsize
	Since there is no control coupling, in order to update the controller for each agent at each iteration we set again,
	\begin{align}
		K_{\text{new}}^{(i)} = -H_{22}^{(i)-1} H_{21}^{(i)}.
		\label{eq:update}
	\end{align}
	Finally, as in the centralized case, for each agent $i$ we define $\phi_t^{(i)} = \bar{z}_t^{(i)} - \gamma\bar{z}_{t+1}^{(i)}$.
	
	%=================================
	
	\subsection{Main Results}
	\label{sec:algorithm}
	
	In this section, we introduce the distributed policy iteration algorithm.
%	The convergence proof and the computational complexity analysis in the next section are inspired by \cite{bradtke1994incremental} that also proposed a similar distributed policy iteration algorithm.
	The analysis in this part is mainly inspired by \cite{bradtke1994incremental}, however, there are fundamental differences as we only assume couplings through a global cost function; as such, the state transition or feedback of each agent only depends on their own history of states and actions.
	Under these assumptions, we show that this way of coupling in the case of identical systems signifies the interdependency of the agents in the decision-making process.
	\begin{algorithm}[H]
		\label{alg:DRL}
		\caption{The distributed Q-learning Algorithm}\label{euclid}
		\begin{algorithmic}[1]
			\State \textbf{Initialize}:
			\State $\hspace{5mm}$ Random: $\hat{\theta}_0^{(1)}(0),\dots,\hat{\theta}_0^{(N)}(0)$
			\State $\hspace{5mm}$ Stabilizable: $K_0^{(1)},\dots,K_0^{(N)}$
			\State $\hspace{5mm}$ $t=0$, $k=1$
			\State \textbf{while convergence:}
			\State $\hspace{5mm}$ Reset Covariance: $P_k(0)=P_0$
			\State $\hspace{5mm}$ \textbf{For} $j=1$ to $M$:
			\State $\hspace{10mm}$  \textbf{For system} $i=1,\dots,N$:
			\State $\hspace{15mm}$ Choose $e_t$ and find $u_t^{(i)}=K_k^{(i)}x_t^{(i)}+e_t$
			\State $\hspace{15mm}$ Collect $x_{t+1}^{(i)}$ by applying $u_t^{(i)}$ to the system
			\State $\hspace{15mm}$ Update $\hat{\theta}_k^{(i)}(j)$ using RLS
			\State $\hspace{15mm}$ $t=t+1$
			\State $\hspace{5mm}$ \textbf{For system} $i=1,\dots,N$:
			\State $\hspace{10mm}$ Find symmetric $\hat{H}_k^{(i)}$ corresponding to $\hat{\theta}_k^{(i)}$
			\State $\hspace{10mm}$ Policy update:  $K_{k+1}^{(i)}=-\hat{H}_{k(22)}^{{(i)}^{-1}}\hat{H}^{(i)}_{k(21)}$
			\State $\hspace{10mm}$ Initialize parameters $\hat{\theta}_{k+1}^{(i)}(0)=\hat{\theta}_k^{(i)}(M)$
			\State $\hspace{5mm}$ $k=k+1$
		\end{algorithmic}
		\label{alg:1}
	\end{algorithm}
	\noindent We briefly explain the steps of the algorithm:
	$\hat{\theta}_k^{(i)}$ is the estimate of $H^{(i)}$ as in \eqref{eq:quad}.
	In the sequel, $\theta^{*(i)}_k$ denotes the parameters of $H^{(i)}$ obtained using the true system parameters.
	$K_k^{(i)}$ denotes the controller estimate.
	The counter $t$ keeps track of the number of collected data while $k$ designates the iteration count on the parameters estimate.
%	Hence by "convergence in the limit" we mean $\lim_{k\rightarrow\infty}$.
	Note that these counters are never reset to zero.
%	The RLS algorithm at each iteration is followed from: \small
%	\begin{align*}
%		\hat{\theta}_k^{(i)}(j) &= \hat{\theta}_k^{(i)}(j-1) + \frac{P_k(j-1)\phi_t^{(i)}\big(r_t^{(i)}-\phi_t^{(i)\top}\hat{\theta}_k^{(i)}(j-1)\big)}{1+\phi_t^{(i)\top}P_k(j-1)\phi_t^{(i)}} \\
%	    P_k(j) &= P_k(j-1) - \frac{P_k(j-1)\phi^{(i)}\phi^{(i)\top}P_k(j-1)}{1+\phi_t^{(i)\top}P_k(j-1)\phi_t^{(i)}} \\
%	    P_k(0) &= P_0
%	\end{align*} \normalsize
%	RLS implementation is obtained from \cite{goodwin2014adaptive}, 
	$P_k(j)$ is the covariance matrix reset to some constant $P_0$ at each iteration to revitalize the gain.
	Each RLS estimation interval includes $M$ time-steps.
	The value of $M$ is dependent on the number of unknown parameters in $\hat{\theta}_k^{(i)}$ and also the desired accuracy.
	The control signal is PE at each iteration of the RLS and $e_t$ is the excitation component which is assumed to be the same for all agents.
	After convergence of RLS, the controller for each agent is updated based on \eqref{eq:update}.
	The estimation parameters are reinitialized from the final value of the previous iteration such that $\hat{\theta}_{k+1}(0) = \hat{\theta}_{k}(M)$.
	The reader is referred to Chapter 3 of \cite{goodwin2014adaptive} for exact steps of RLS.
%	The convergence proof of the similar algorithm for the centralized case is given in \cite{bradtke1994incremental}
%	In the next theorem, we will show that for the case of identical dynamically decoupled interconnected systems, for each agent the algorithm leads to the corresponding optimal controller.
	
	\vspace{2mm}
	
	\begin{theorem}
		Assume that for all $i=1,\dots,N$, the pair  $(A,B)$ is a controllable and $K_0^{(i)}$ is stabilizing with a PE signal $\phi_t^{(i)}$.
		Then there exists $M<\infty$ such that \Cref{alg:1} generates a sequence $\{ K_k^{(i)}\}$ with $\lim_{k\rightarrow\infty} \|K_k^{(i)} - K^{*}\| = 0$, where $K^{*}=\text{LQR}(A,B,\bar{Q},\bar{R})$.
	\end{theorem}
	
	\vspace{2mm}
	
	\begin{proof}
		From \eqref{eq:reward},
		\begin{align*}
			r^{(i)}_t = y_t^{(i)\top} H^{(i)} y_t^{(i)} - \gamma y_{t+1}^{(i)\top} H^{(i)} y_{t+1}^{(i)}.
		\end{align*}
		Also from \cref{sec:DQfun},
		\begin{align*}
			r_t^{(i)} &= x_t^{(i)\top} \bar{Q} x_t^{(i)} + u_t^{(i)\top} \bar{R} u_t^{(i)} \\
			&\hspace{15mm} + \sum_{k=1}^{d_i} \big( x_t^{(i)}-x_t^{(j_k)} \big)^{\top} \bar{Q} \big( x_t^{(i)}-x_t^{(j_k)} \big),
		\end{align*}
		and, \footnotesize
		\begin{equation*}
			\begin{aligned}
				y_t^{(i)\top} & H^{(i)} y_t^{(i)} \\
				&= \big[x_t^{(i)\top}\ u_t^{(i)\top} \big]
				\begin{bmatrix}
					\bar{Q}+\gamma A^{\top}P^{(i)}A & \gamma A^{\top}P^{(i)}B \\
					\gamma B^{\top}P^{(i)}A & \bar{R}+\gamma B^{\top}P^{(i)}B 
				\end{bmatrix}
				\begin{bmatrix}
					x_t^{(i)} \\ u_t^{(i)}
				\end{bmatrix} .\\
%				& \hspace{28mm} + \sum_{k=1}^{d_i} \big( x_t^{(i)}-x_t^{(j_k)} \big)^{\top} \bar{Q} \big( x_t^{(i)}-x_t^{(j_k)} \big) \\
%				y_{t+1}^{(i)\top} & H^{(i)} y_{t+1}^{(i)} \\
%				&= \big[x_{t+1}^{(i)\top}\ u_{t+1}^{(i)\top} \big]
%				\begin{bmatrix}
%					\bar{Q}+\gamma A^{\top}P^{(i)}A & \gamma A^{\top}P^{(i)}B \\
%					\gamma B^{\top}P^{(i)}A & \bar{R}+\gamma B^{\top}P^{(i)}B 
%				\end{bmatrix}
%				\begin{bmatrix}
%					x_{t+1}^{(i)} \\ u_{t+1}^{(i)}
%				\end{bmatrix} \\
%				& \hspace{27mm} + \sum_{k=1}^{d_i} \big( x_{t+1}^{(i)}-x_{t+1}^{(j)} \big)^{\top} \bar{Q} \big( x_{t+1}^{(i)}-x_{t+1}^{(j)} \big).
			\end{aligned}
		\end{equation*} \normalsize
		Consequently, we obtain,
		\begin{equation}
			\begin{aligned}
				\xi_t^{(i)} = \bar{z}_t^{(i)\top} \hat{\theta}^{(i)}_k
			\end{aligned}
			\label{eq:RLSagent}
		\end{equation}
		where,
		\begin{equation}
			\begin{aligned}
				\xi_t^{(i)} = x_t^{(i)\top}& \bar{Q} x_t^{(i)} + u_t^{(i)\top} \bar{R} u_t^{(i)} \\
				&+ \sum_{k=1}^{d_i} \big( x_{t+1}^{(i)}-x_{t+1}^{(j)} \big)^{\top} \bar{Q} \big( x_{t+1}^{(i)}-x_{t+1}^{(j)} \big).
			\end{aligned}
			\label{eq:introCoupling}
		\end{equation}
		Hence the distributed nature of the problem narrows down to a particular distributed form of RLS.
%		We will show that after some iteration, the terms $x_{t+1}^{(i)}-x_{t+1}^{(j)}$ converge to zero as a result of the convergence of the controller estimates.
		We stack the $N$ equations of the form \eqref{eq:RLSagent} for all agents into vector form as,
		\begin{align}
		\underbrace{\begin{bmatrix}
			\xi^{(1)}_t \\ \vdots \\ \xi^{(N)}_t
			\end{bmatrix}}_{\Xi_t} =
		\underbrace{\begin{bmatrix}
			\bar{z}_t^{(1)^{\top}} & & \\
			& \ddots & \\
			& & \bar{z}_t^{(N)^{\top}}
			\end{bmatrix}}_{Z_t}
		\underbrace{\begin{bmatrix}
			\hat{\theta}_k^{(1)} \\ \vdots \\ \hat{\theta}_k^{(N)}
			\end{bmatrix}}_{\hat{\Theta}_k} .
		\label{eq:RLS}
		\end{align}
%		where,
%		\begin{equation}
%			\begin{aligned}
%				\Xi_t &= \begin{bmatrix}
%					x_t^{(1)\top} Q x_t^{(1)} + u_t^{(1)\top} R u_t^{(1)} \\ \vdots \\ x_t^{(N)\top} Q x_t^{(N)} + u_t^{(N)\top} R u_t^{(N)}
%				\end{bmatrix} \\
%				&\hspace{16mm} + \begin{bmatrix}
%					\sum_{k=1}^{d_1} \big( x_{t+1}^{(1)}-x_{t+1}^{(k)} \big)^{\top} Q \big( x_{t+1}^{(1)}-x_{t+1}^{(k)} \big)
%					\\ \vdots \\
%					\sum_{k=1}^{d_N} \big( x_{t+1}^{(N)}-x_{t+1}^{(k)} \big)^{\top} Q \big( x_{t+1}^{(N)}-x_{t+1}^{(k)} \big)
%				\end{bmatrix}
%			\end{aligned}
%			\label{eq:bigRLS}
%		\end{equation} \normalsize
%		We will show that for identical systems, the second term in \eqref{eq:bigRLS} will go to zero.
%		This is not trivial since the output of the RLS includes data that comes from different agents in the network.
		Based on the definition of PE in \eqref{eq:PE}, it is straightforward to show that the matrix $Z_t$ is PE if $\bar{z}_t^{(i)}$ is PE for all $i$.
		This results in the convergence of \cref{eq:RLS} to some $\Theta^*$ for large enough $M$.\footnote{Parameter estimation for the multi-output system is an straightforward extension of the scalar case and is discussed in Chapter 3.8 of \cite{goodwin2014adaptive}.}
		From Theorem 5.1 in \cite{bradtke1994incremental},
		\begin{equation}
			\begin{aligned}
				\lim_{k\rightarrow\infty} \| \hat{\theta}^{(i)}_k - \theta^{*(i)}_k \| = 0, \quad
				\lim_{k\rightarrow\infty} \| \theta^{*(i)}_k - \theta^{*(i)}_{k-1} \| = 0.
			\end{aligned}
			\label{eq:conv0}
		\end{equation}
		However, the convergence of $\hat{\theta}_k^{(i)}$ for all $i$ to one single value is non-trivial due to the interdependency in RLS.
		We will show that for a connected network of agents,
		\begin{align*}
			\lim_{t\rightarrow \infty}\ \| x_t^{(i)} - x_t^{(j)} \| = 0,
		\end{align*}
		for any $i$ and $j$.
		Note that according to \eqref{eq:introCoupling}, if a node is disconnected from the graph it can be individually examined as in the centralized case.
		%Also the limit in $t$ implies the increase in $k$.
		Recall that for $\ell = i, j$,
		\begin{align*}
			x_{t+1}^{(\ell)} = Ax_t^{(\ell)} + Bu_t^{(\ell)} = (A-BK_k^{(\ell)})x_t^{(\ell)} + Be_t.
		\end{align*}
		As such,
		\begin{align}
			x_{t+1}^{(i)} - x_{t+1}^{(j)} = (A-BK_k^{(i)}) (x_t^{(i)} - x_t^{(j)}) + B \Delta K_k x_t^{(j)}
			\label{eq:difference}
		\end{align}
		where $\Delta K_k = K_k^{(i)} - K_k^{(j)}$.
		Then if we show that $\|\Delta K_k\|\rightarrow 0$ as $k\rightarrow \infty$ we obtain,
		\begin{align}
			\| x_{t+1}^{(i)} - x_{t+1}^{(j)} \| = \| (A-BK_k^{(i)})^{t+1} (x_0^{(i)} - x_0^{(j)}) \| \rightarrow 0,
			\label{eq:final}
		\end{align}
		given that the policy iteration algorithm leads to a more stabilizing controller $K_k^{(i)}$ as $k$ increases \cite{bradtke1994incremental}.
		Then,
		\begin{equation}
			\begin{aligned}
				\| K_k^{(i)} &- K_k^{(j)} \| \\
				&= \| \hat{H}_{{k-1}_{(22)}}^{{(j)}^{-1}}\hat{H}^{(j)}_{{k-1}_{(21)}} - \hat{H}_{{k-1}_{(22)}}^{{(i)}^{-1}}\hat{H}^{(i)}_{{k-1}_{(21)}} \| \\
%				&= \| \hat{H}_{{k-1}_{(22)}}^{{(j)}^{-1}}\hat{H}^{(j)}_{{k-1}_{(21)}} - \hat{H}_{{k-1}_{(22)}}^{{(j)}^{-1}}\hat{H}^{(i)}_{{k-1}_{(21)}} \\
%				&\hspace{10mm}+ \hat{H}_{{k-1}_{(22)}}^{{(j)}^{-1}}\hat{H}^{(i)}_{{k-1}_{(21)}} - \hat{H}_{{k-1}_{(22)}}^{{(i)}^{-1}}\hat{H}^{(i)}_{{k-1}_{(21)}} \| \\
				&= \| \hat{H}_{{k-1}_{(22)}}^{{(j)}^{-1}} \Big( (\hat{H}^{(j)}_{{k-1}_{(21)}} - \hat{H}^{(i)}_{{k-1}_{(21)}})\\
				&\hspace{10mm} + (\hat{H}^{(i)}_{{k-1}_{(22)}} - \hat{H}^{(j)}_{{k-1}_{(22)}}) \hat{H}_{{k-1}_{(22)}}^{{(i)}^{-1}} \hat{H}^{(i)}_{{k-1}_{(21)}} \Big) \| .
			\end{aligned}
			\label{eq:HH}
		\end{equation}
		Since $\hat{H}_{22}$ and $\hat{H}_{21}$ contain only a subset of elements in $\hat{\theta}$,
		\begin{equation}
			\begin{aligned}
				\| \hat{H}^{(j)}_{{k-1}_{(21)}} - \hat{H}^{(i)}_{{k-1}_{(21)}} \| &\leq \| \hat{\theta}_{k-1}^{(j)} - \hat{\theta}_{k-1}^{(i)} \|, \\ 
				\| \hat{H}^{(j)}_{{k-1}_{(22)}} - \hat{H}^{(i)}_{{k-1}_{(22)}} \| &\leq \| \hat{\theta}_{k-1}^{(j)} - \hat{\theta}_{k-1}^{(i)} \|.
			\end{aligned}			
			\label{eq:ineq}
		\end{equation}
		Hence equations \eqref{eq:HH} and \eqref{eq:ineq} lead to, 
		\begin{equation}
			\begin{aligned}
				\| & K_k^{(i)} - K_k^{(j)} \| \\
				&\leq \| \hat{H}_{{k-1}_{(22)}}^{{(j)}^{-1}} \|\ \| \hat{\theta}_{k-1}^{(j)} - \hat{\theta}_{k-1}^{(i)} \|\ \| 1 + \hat{H}_{{k-1}_{(22)}}^{{(i)}^{-1}} \hat{H}^{(i)}_{{k-1}_{(21)}} \| \\
				&\leq \kappa_0 \| \hat{\theta}_{k-1}^{(j)} - \hat{\theta}_{k-1}^{(i)} \| ,
			\end{aligned}
			\label{eq:uBound}	
		\end{equation}
		where we have used the fact that the estimated parameters are bounded and $\kappa_0 > 0$ is a constant such that,
		\begin{align*}
			\| \hat{H}_{{k-1}_{(22)}}^{{(j)}^{-1}} \| ~.~ \| 1 + \hat{H}_{{k-1}_{(22)}}^{{(i)}^{-1}} \hat{H}^{(i)}_{{k-1}_{(21)}} \| \leq \kappa_0.
		\end{align*}
		From Lemma 5.2 in \cite{bradtke1994incremental},
		\begin{align*}
			\| \theta_k^{*(\ell)} - \hat{\theta}_k^{(\ell)} \| \leq \epsilon_M \big( \| \theta_k^{*(\ell)} -  \theta_{k-1}^{*(\ell)}\| + \| \theta_{k-1}^{*(\ell)} - \hat{\theta}_{k-1}^{(\ell)} \| \big),
		\end{align*}
		which for large enough $M$ results in,
		\begin{equation}
		\begin{aligned}
			\label{eq:conv1}
			&\| \theta_k^{*(i)} - \hat{\theta}_k^{(i)} \| + \| \theta_k^{*(j)} - \hat{\theta}_k^{(j)} \| \\
			&\hspace{13mm} \leq \epsilon_M \big( \| \theta_k^{*(i)} -  \theta_{k-1}^{*(i)}\| + \| \theta_{k-1}^{*(i)} - \hat{\theta}_{k-1}^{(i)} \| \\
			&\hspace{19.5mm} + \| \theta_k^{*(j)} -  \theta_{k-1}^{*(j)}\| + \| \theta_{k-1}^{*(j)} - \hat{\theta}_{k-1}^{(j)} \| \big).
		\end{aligned}
		\end{equation}
		Using triangle inequality on the left side of this inequality,
		\begin{equation}
		\begin{aligned}
			\label{eq:conv2}
			&\Big| \| \hat{\theta}_k^{(i)} - \hat{\theta}_k^{(j)} \| - \| \theta_k^{*(i)} - \theta_k^{*(j)} \| \Big| \\
			&\hspace{23mm}\leq \big\| \big( \hat{\theta}_k^{(i)} - \hat{\theta}_k^{(j)} \big) - \big( \theta_k^{*(i)} - \theta_k^{*(j)} \big) \big \| \\
			&\hspace{23mm}\leq \| \theta_k^{*(i)} - \hat{\theta}_k^{(i)} \|\ +\ \| \theta_k^{*(j)} - \hat{\theta}_k^{(j)} \|.
		\end{aligned}
		\end{equation}
		Then, from \eqref{eq:conv1} and \eqref{eq:conv2} and for large $k$,
		\begin{align*}
			\| \hat{\theta}_k^{(i)} - \hat{\theta}_k^{(j)} \| &\leq \epsilon_M \Big( \| \theta_k^{*(i)} -  \theta_{k-1}^{*(i)}\| + \| \theta_{k-1}^{*(i)} - \hat{\theta}_{k-1}^{(i)} \| \\
			&\hspace{7mm} +\| \theta_k^{*(j)} -  \theta_{k-1}^{*(j)}\| + \| \theta_{k-1}^{*(j)} - \hat{\theta}_{k-1}^{(j)} \| \Big) \\
			&\hspace{7mm} + \| \theta_k^{*(i)} - \theta_k^{*(j)} \|.
		\end{align*}
		Hence using the result in \eqref{eq:conv0},
		\begin{align*}
			\| \hat{\theta}_k^{(i)} - \hat{\theta}_k^{(j)} \| \rightarrow 0,
		\end{align*}
		and plugging this into \eqref{eq:uBound},
		\begin{align}
			\|\Delta K_k\| = \| K_k^{(i)} - K_k^{(j)} \| \rightarrow 0,
		\end{align}
		Hence,
		\begin{align*}
		\| x_{t+1}^{(i)} - x_{t+1}^{(j)} \| \rightarrow 0.
		\end{align*}
		This implies that based on \eqref{eq:introCoupling}, for identical systems the algorithm moves towards $N$ decoupled $Q$-learning algorithms for each agent.
		Thus, although the provided data is from an interconnected system, each controller converges to its optimal value, i.e.,
		\begin{align*}
			\lim_{k\rightarrow\infty} \|K_k^{(i)} - K^{*}\| = 0, \qquad \text{for} \quad i=1,2,\dots,N .
		\end{align*}
	\end{proof}
	
	\begin{remark}
		In \Cref{alg:1}, we have assumed that the exploration signal, $e_t$, is equal for every agent at each time step.
		This is a valid assumption as long as $Z_t$ in \eqref{eq:RLS} is PE so that RLS is assured to converge.
		Another option would be to choose the excitation signals $e_t^{(i)}$ and $e_t^{(j)}$ in a way that,
		\begin{align*}
			e_t^{(i)} - e_t^{(j)} = -\Delta K_k x_t^{(j)}.
		\end{align*}
		Hence, not only the input to the RLS is PE, the difference cancels out $\Delta K_k x_t^{(j)}$ in \eqref{eq:difference}.
		However, this setup is more challenging to implement, particularly for large-scale systems.
	\end{remark}
	
	%=================================
	
	\subsection{Computational Saving}
	\label{sec:comSaving}
	
	The computational saving resulting from using the distributed $Q$-learning algorithm is significant, since for a large system, the design of the LQR controller with the computational complexity of solving ARE of order $\mathcal{O}(n^3)$, can be prohibitively expensive.
	The main computational burden of \Cref{alg:1} comes from RLS where the complexity of its implementation is $\mathcal{O}(\gamma^2)$ with $\gamma$ parameters to learn.
	Assuming that the system contains $N$ agents each having $n$ states and $m$ inputs, the computational complexity of the centralized $Q$-learning is obtained by,
	\begin{align*}
		\mathcal{O}\Big( \big( \frac{(Nn+Nm)(Nn+Nm+1)}{2} \big)^2 \Big) \approx \mathcal{O} \big( N^4(n+m)^4 \big),
	\end{align*}
	while for the distributed case the code performs $N$ repetitions of the same RLS leading to the complexity bound,
	\begin{align*}
		\mathcal{O}\Big( N \big( \frac{(n+m)(n+m+1)}{2} \big)^2 \Big) \approx \mathcal{O} \big( N(n+m)^4 \big).
	\end{align*}
	Hence the complexity reduction is,
	\begin{align*}
		\frac{N^4(n+m)^4 - N(n+m)^4}{N^4(n+m)^4} \times 100 = \frac{N^3-1}{N^3} \times 100 \% ,
	\end{align*}
which is substantial for large $N$.
	\Cref{tab:save} compares the computational saving for some values of $N$.
	\begin{table}[h]
		\centering
		\resizebox{0.49\textwidth}{!}{
			\begin{tabular}{|c||c|c|c|c|c|}
				\hline
				N & 2 & 3 & 5 & 8 & 100 \\ [0.5ex] 
				\hline
				Saving (\%) & 87.5 & 96.29 & 99.2 & 99.8 & 99.99 \\ [0.5ex]
				\hline
			\end{tabular}
		}
		\caption{Approximate computational saving of the distributed $Q$-learning compared to the centralized case in \cite{bradtke1994incremental}.}
		\label{tab:save}
	\end{table}
	
	%=================================
	
	\section{Example}
	\label{sec:example}
	
\begin{figure}[H]
		\centering
		\includegraphics[scale=0.19]{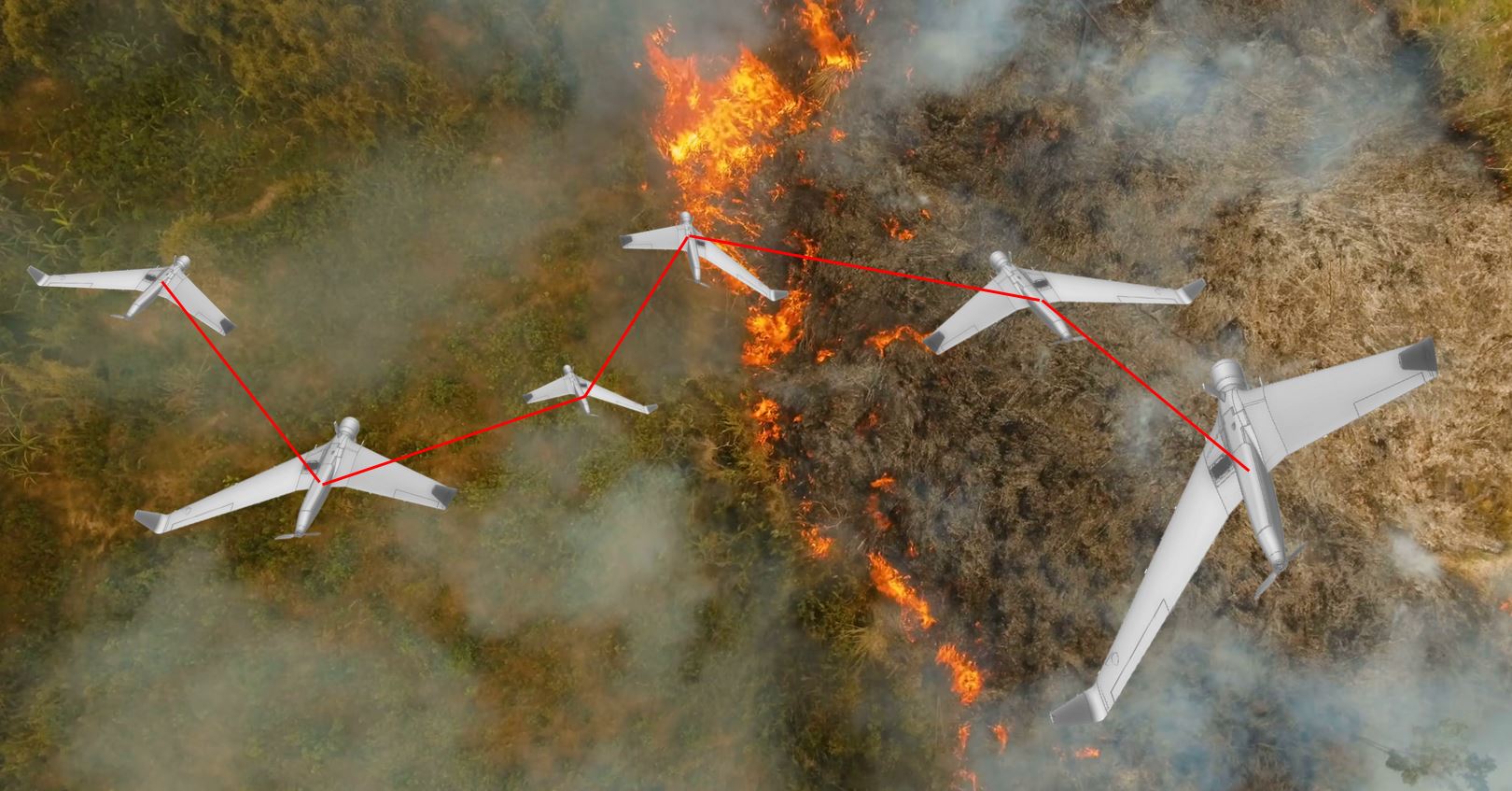}
		\caption{A group of identical firefighting UAVs maneuvering in parallel aiming to extinguish a blaze (Aerial view of the forest fire - \emph{Photo Credit: Alex Punker, Bigstock}).}
		\label{fig:1}
	\end{figure}

\begin{figure}[H]
		\centering
		\includegraphics[scale=0.16]{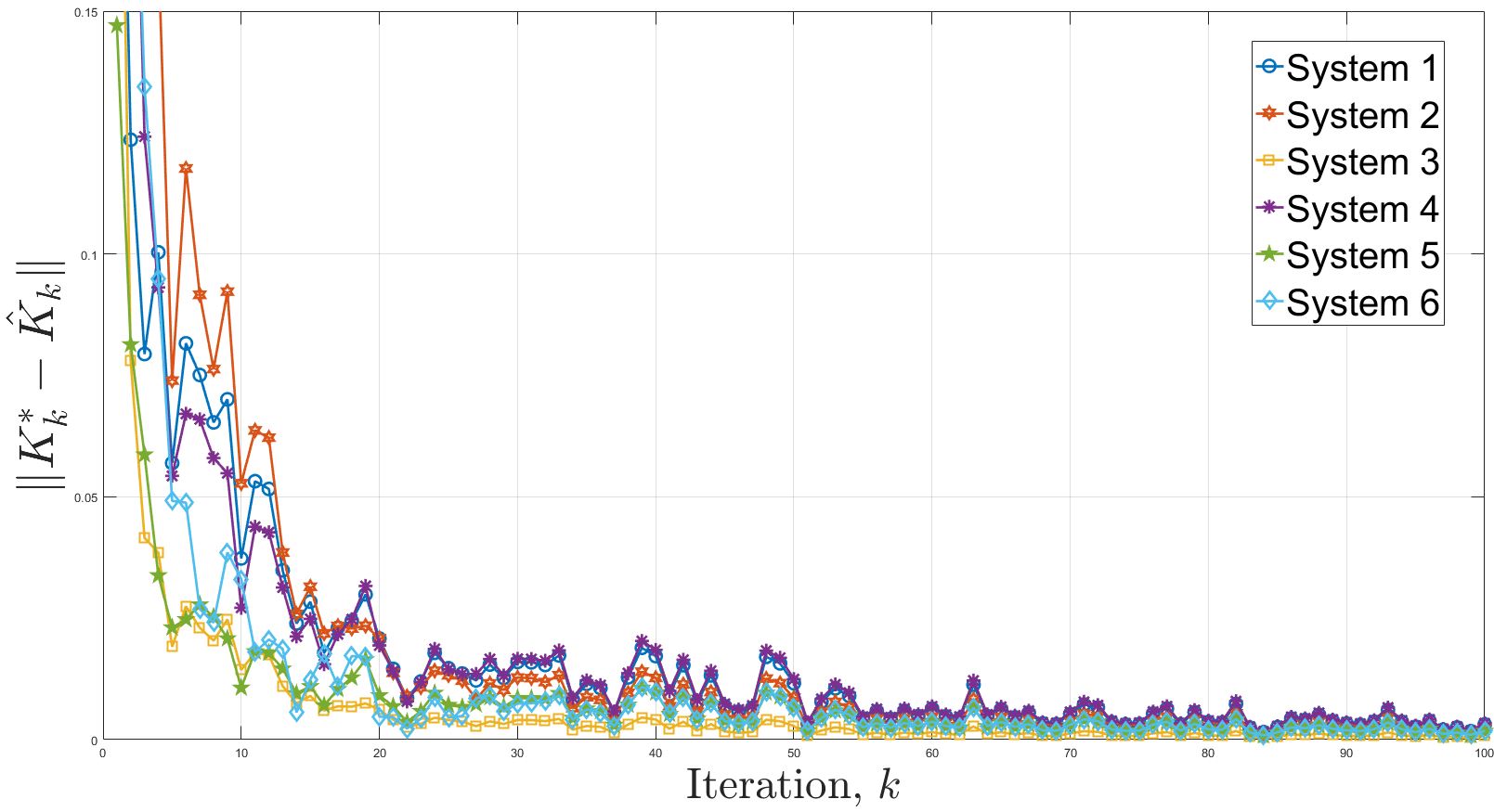}
		\caption{Performance of the distributed Q-learning algorithm. 
		The plot demonstrates the norm of the error between the LQR optimal controller of each subsystem and the estimate of the algorithm at each iteration $k$.}
		\label{fig:2}
	\end{figure}

\begin{figure}[H]
		\centering
		\includegraphics[scale=0.16]{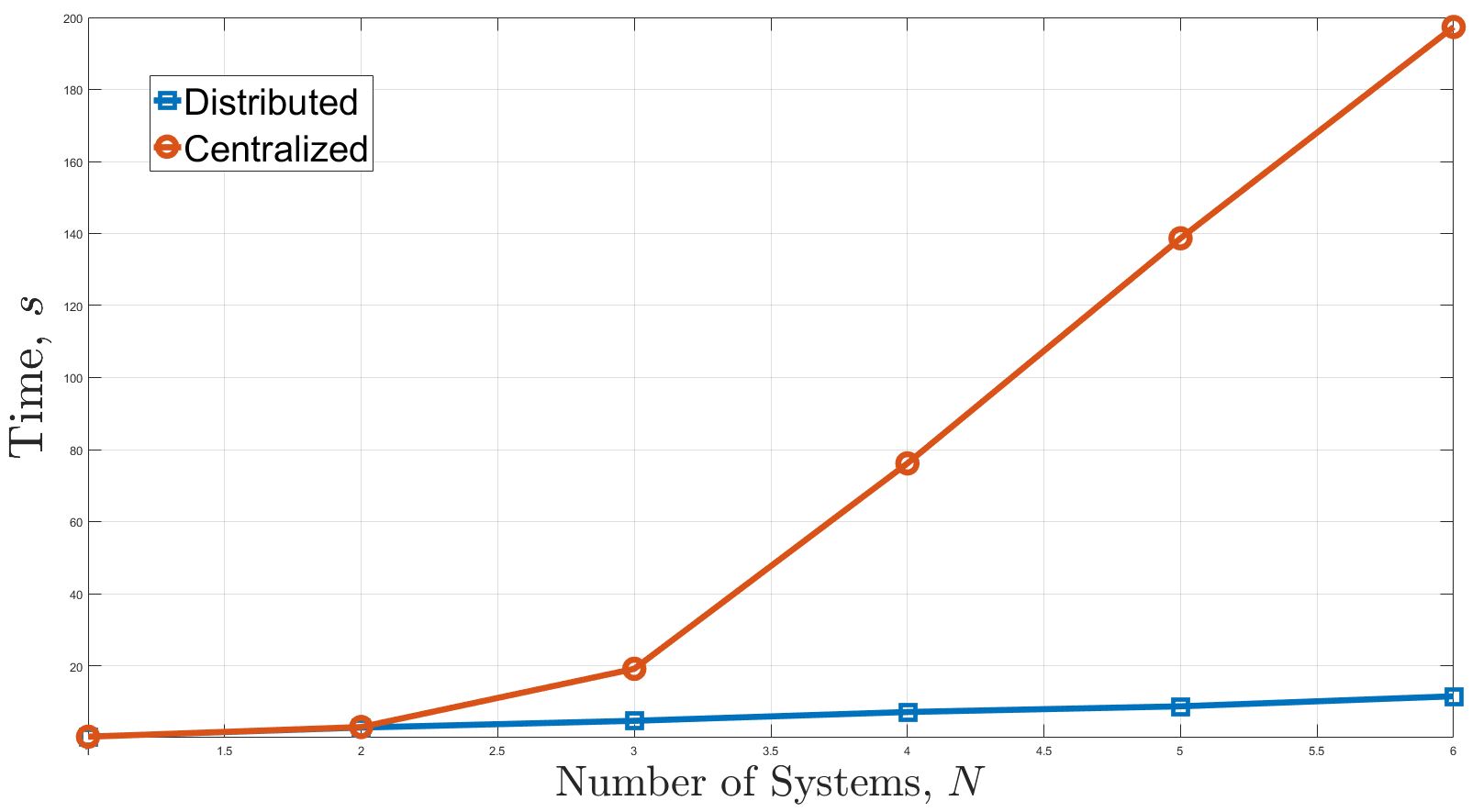}
		\caption{Computational performance of distributed and centralized algorithm for different number of interacting agents.}
		\label{fig:3}
	\end{figure}

	In this section, we provide an example to show the efficiency of the distributed $Q$-learning algorithm for a set of identical communicating UAVs.
	We consider the autonomous flight of a network of six interconnected Unmanned Aerial Vehicles (UAVs) which are set to perform a common task such as geographical data collection or putting out a wildfire.
	To cover the whole targeted area, these UAVs are programmed to move in parallel and in order for minimal signal transmissions, each UAV only communicates with its closest neighbor in the network as depicted in \Cref{fig:1}.
	\noindent The discrete-time dynamics of UAVs is considered by mini-aircraft linear parameters that can be found in \cite{hung1982multivariable}.
%	the dynamics of each UAV is defined by,
%	\begin{align*}
%		A &= \begin{bmatrix}
%			0.7000  &  0.0014  &  0.1132  &  0.0005 &  -0.0967 \\
%			0   &   0.6945 &  -0.0171  & -0.0005  &  0.0068 \\
%			0   &   0.0003  &  0.7000  &  0.0957 &  -0.0048 \\
%			0   &   0.0060 &  -0.0000  &  0.6131  & -0.0936 \\
%			0   &  -0.0277  &  0.0002  & 0.0973  &  0.6287
%		\end{bmatrix} \\
%		B &= \begin{bmatrix}
%		   -0.0076  &  0.0000  &  0.0003 \\
%   		   -0.0115  &  0.0997  &  0.0000 \\
%			0.0212  &  0.0000  & -0.0081 \\
%			0.4152  &  0.0003  & -0.1589 \\
%			0.1742  & -0.0014  & -0.0154
%		\end{bmatrix}
%	\end{align*}
%	for a discrete-time linear time-invariant system.
%	The states and inputs represent:
%	\begin{itemize}
%		\item $x_1$: altitude (m)
%		\item $x_2$: forward speed (m/s)
%		\item $x_3$: pitch angle (deg)
%		\item $x_4$: pitch rate(deg/s)
%		\item $x_5$: vertical speed (m/s)
%		\item $u_1$: spoiler angle (deg/10)
%		\item $u_2$: forward acceleration (m/$s^2$)
%		\item $u_3$: elevator angle (deg)
%	\end{itemize}
	We assume $\bar{Q}=\textsc{I}_5$ and $\bar{R} = \textsc{I}_3$.
	We will show the results of the distributed policy iteration for $N=6$, $n=5$, and $m=3$ and compare the computational performance with the centralized case.
	For the distributed algorithm we consider $M=50$ and the exploration signal $e_t$ is generated from a normal distribution.
	\Cref{fig:2} shows the results of simulations regarding the controller error norm.
	\noindent A comparison between the computational performance of the centralized and distributed methods is also provided in \Cref{fig:3}.
	For scaling purposes, $M$ and $e_t$ are re-adjusted for each $N$.\footnote{The main reason for this is that $M$ needs to be modified since $N$ is proportionally related to the centralized system dimensions $Nn$ and $Nm$.}
%	The results verify the computational saving which was discussed in \Cref{sec:comSaving}.
	
	%=================================
	
	\section{Conclusion}
	\label{sec:conclusion}
	
	In this paper we examined a data-guided approach for the control of large-scale interconnected identical systems with decoupled dynamics; it is assumed that the interconnection is reflected in the cost function for the control design problem.
	We leveraged a distributed $Q$-learning as a policy iteration method. In this direction, it is shown that the proposed distributed algorithm converges to each agent's individual optimal controller,  which could have been obtained by running a centralized $Q$-learning algorithm.
	The significance of the resulting computational savings are also discussed.
	
	There are a number of directions to pursue as future works.
	First, the observation in this paper can be further extended to more elaborate cost structure, highlighting the trade-off between local and global optimality in large-scale distributed systems.
	This can be achieved if other types of interconnections such as dynamics or feedback coupling as well as consensus through the $Q$-function are adopted for the analysis.\footnote{A similar approach has been examined in \cite{kar2013cal}.}
	Another line of work is to consider other types of data-guided distributed control mechanisms for structures such as layering or systems with switching dynamics.
		%=================================
	
	\bibliographystyle{ieeetr}
	\bibliography{citations2}
	
\end{document}